\documentclass{article}

\usepackage{amssymb}
\usepackage{psfig}
\usepackage{amscd}

\newcommand{\rr}{{\mathbb R}}
\newcommand{\cc}{{\mathbb C}}
\newcommand{\bu}{{\bf u}}
\newcommand{\bx}{{\bf x}}
\newcommand{\by}{{\bf y}}
\newcommand{\bz}{{\bf z}}

\begin{document}

\title{Computing Dynamic Output Feedback Laws}

\author{Jan Verschelde\thanks{
Department of Mathematics, Statistics, and Computer Science,
University of Illinois at Chicago, 851 South Morgan (M/C 249),
Chicago, IL 60607-7045, USA.
{\em Email:} jan@math.uic.edu or jan.verschelde@na-net.ornl.gov.
{\em URL:} http://www.math.uic.edu/{\~{}}jan.
This material is based upon work
supported by the National Science Foundation under Grant No.\
0105739 and Grant No.\ 0134611.}
\and
Yusong Wang\thanks{
Department of Mathematics, Statistics, and Computer Science,
University of Illinois at Chicago, 851 South Morgan (M/C 249),
Chicago, IL 60607-7045, USA.
{\em Email:} ywang25@uic.edu.}}

\date{2 June 2003}

\maketitle

\begin{abstract} 
\noindent The pole placement problem asks to find laws to feed the
output of a plant governed by a linear system of differential equations
back to the input of the plant so that the resulting closed-loop 
system has a desired set of eigenvalues.
Converting this problem into a question of enumerative geometry,
efficient numerical homotopy algorithms to solve this problem
for general Multi-Input-Multi-Output (MIMO) systems have been
proposed recently.
While dynamic feedback laws offer a wider range of use,
the realization of the output of the numerical homotopies as
a machine to control the plant in the time domain has not been
addressed before.  In this paper we present symbolic-numeric
algorithms to turn the solution to the question of enumerative
geometry into a useful control feedback machine.
We report on numerical experiments with our publicly available
software and illustrate its application on various control problems
from the literature.



\end{abstract}

\section{Introduction}

Given a linear system of differential equations
and a list of eigenvalues,
the pole placement problem asks to find laws to feed the
output back to the input so that the resulting closed-loop
system has the same eigenvalues as the given list.
While executing a static feedback law corresponds to a 
simple matrix-vector multiplication, a dynamic compensator            
has several internal states.

The theoretical solution of this problem was found in 
the Schubert calculus, see~\cite{BB81,Byr89} for the static,
and~\cite{RRW96,RRW98} for dynamic compensators.
Because of its importance to practical applications,
the development of algorithms for this problem was stated
as an open problem~\cite{RW98} (see also~\cite{Chu01}).
The first homotopy algorithms were proposed in~\cite{HSS98},
defining a numerical Schubert calculus.
The Pieri homotopy algorithms of~\cite{HSS98} were improved and
generalized to dynamic feedback in~\cite{HV00}.  In~\cite{LWW02}, 
the numerical performance of these homotopies for static feedback
was improved.

Our current paper is a sequel to our previous paper~\cite{VW02}
as a continuing effort to apply the homotopy algorithms to
practical applications of linear systems control.  Rather than
generating random matrices as inputs to our algorithms we look
to build case studies with examples from the literature.
In~\cite{VW02} we applied static output feedback to use pole
placement to keep a satellite in orbit.  In this current paper
we consider the application of dynamic feedback laws to this
and other applications.  See~\cite{KRW03} for the relation of
the pole placement problem with inverse eigenvalue and 
matrix extension problem.

As expected, with the more general dynamic feedback laws we
can cover a wider range of applications, but at the same time
we experienced the real need for more advanced algorithms from
numerical algebraic geometry to deal better with over and 
underdetermined problems.  We define these mathematical 
problems in the next section.

The heart of this paper concerns the realization of dynamic
compensators.  As these compensators are computed in the
frequency domain, at the end of the paths defined by the
homotopies, any application demands the description of the
controllers in the time domain as a tuple of matrices.
While realization algorithms are covered in standard graduate textbooks
(e.g.,~\cite{AN97} and~\cite{Kai80}), our output is numerical, 
i.e.: subject to roundoff and known only with limited accuracy.
Reports on numerical implementations of algorithms which manipulate
polynomial matrices are scarce, we found only one paper~\cite{HS99},
which presents a numerical study on a commercial implementation,
but leaves the
implementation of a numerical Smith normal form as an open question.
Note that the Smith normal form plays an important role 
to obtain a minimal realization of the 
transfer function of Multi-Input-Multi-Output (MIMO) systems.

Section three starts with the description of an algorithm to
compute the greatest common divisor of two polynomials with
approximate coefficients, as this is the cornerstone function
in the implementation of a numerical Smith normal form.
As we experienced and illustrate, the naive application of
the Euclidean algorithm for an approximate GCD can fail miserably.
That many algorithms to solve problems with exact data turn out
to be numerically unstable is a growing concern in computer algebra,
which has led to hybrid symbolic-numeric computation~\cite{CKW02}.
In particular, the approximate GCD problem has received a lot
of research attention in recent years, see for instance~\cite{BL00},
\cite{EGL97}, \cite{HKL99}, \cite{KL96}, \cite{Pan01},
\cite{ZN00} and~\cite{ZW98}.

Our approach is to consider the input polynomials not as polynomials
with approximate coefficients, but as polynomials defined by 
approximate
roots.  So we view polynomials not symbolically (as sums of monomials
with approximate coefficients) but geometrically (as the collection
of roots in the complex plane, known with limited accuracy).
The idea to computing the GCD by matching common approximate roots
(within a certain tolerance) can be found in~\cite{Pan01}.
Compared (see~\cite{Pan01}) to alternative constructions of approximate 
GCDs,
this method works well for input polynomial with multiple or clustered
zeros.  Algorithmically, we can immediately apply any numerical root
finding method (see~\cite{Pan97} for a recent survey).
To compute the polynomials in the extended GCD problem 
(i.e.: the $k(x)$ and $l(x)$ in the combination 
$k(x) a(x) + l(x) b(x) = {\rm GCD}(a(x),b(x))$), 
we apply Newton interpolation.

As usual,
the link between algorithms and applications is software.
We describe an extension to the software package PHCpack~\cite{Ver99}
in section five before illustrating its application on various
control problems from the literature.

\section{Problem Statement}

The algorithmic framework is defined by the transitions between
the time and the frequency domain (which leads to the symbolic-numeric
calculations).  Mathematically, we count the number of conditions
to determine the dimension of our solution space.

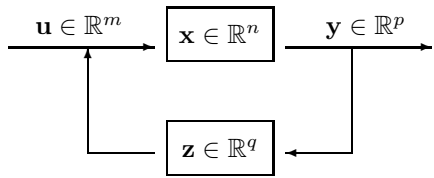
\begin{figure}[hbt]
\begin{center}
\begin{picture}(200,70)(0,0)

\put(30,60){${\bf u} \in \rr^m$}
\put(20,55){\vector(1,0){55}}
\put(80,50){\framebox(40,20)[c]{${\bf x} \in \rr^n$}}
\put(125,55){\vector(1,0){55}}
\put(140,60){${\bf y} \in \rr^p$}

\put(150,15){\vector(-1,0){25}}
\put(150,15){\line(0,1){40}}
\put(80,7){\framebox(40,20)[c]{${\bf z} \in \rr^q$}}
\put(50,15){\line(1,0){25}}
\put(50,15){\vector(0,1){40}}

\end{picture}
\caption{Control of an $m$-input and $p$-output plant by a $q$th order
         dynamic compensator in the time domain.}
\label{figmachine}
\end{center}
\end{figure}

In the time domain, we picture the control of a plant with a dynamic
feedback law in Figure~\ref{figmachine}.

The machine we want to control is given by three matrices $(A,B,C)$,
defining a system of first-order differential equations in the
time domain.  The input of the Pieri homotopies is sampled in the
frequency domain.  So to apply these homotopies, we need to
translate the input from the time to frequency domain and 
to realize the output as a tuple of matrices.
Table~\ref{tabdomains} defines the three stages in the data flow.
The labels on the arrows in Figure~\ref{figdomains}
correspond to the numbers in Table~\ref{tabdomains}.

\begin{table}[hbt]
\begin{center}
\begin{tabular}{|ll|} \hline
  (1) & With plain linear algebra we compute the input for the Pieri \\
      & homotopies, sampling points from the plant we wish to control. \\
      \hline
  (2) & Pieri homotopies compute solution maps of degree~$q$, \\
      & which are then brought into another format. \\ \hline
  (3) & Given the transfer function of the compensator, we realize  \\
      & the compensator by a tuple of four matrices $(F,G,H,K)$. \\ 
\hline
\end{tabular}
\caption{Transitions in applying the Pieri homotopies.}
\label{tabdomains}
\end{center}
\end{table}

Our problem is thus to process the output of the Pieri homotopy
algorithms and to apply realization algorithms to the numerical output.

In~\cite{HV00}, the following equivalence relation was derived:
\begin{eqnarray} \label{eqrewrite}
 \lefteqn{
  \det\left(
    s \left[
         \begin{array}{cc}
            I_n & 0 \\ 0 & I_q
         \end{array}
      \right]
    - \left[
         \begin{array}{cc}
            A + BKC & BH \\ GC & F
         \end{array}
      \right]
    \right) = 0} \\
  & \Leftrightarrow  &
  \det\left[
         \begin{array}{cc}
             I_p & C ( s I_n - A )^{-1} B \\
             H ( s I_q - F )^{-1} G + K & I_m
         \end{array}
      \right] = 0  \label{eqgeometric}
\end{eqnarray}
Equation~(\ref{eqrewrite}) is the characteristic equation of
the closed-loop system.  Via elementary row and column operations,
this equation can be rewritten into~(\ref{eqgeometric}), which
separates the given data $(A,B,C)$ from the unknown $(F,G,H,K)$.
Equation~(\ref{eqgeometric}) shows the geometric problem:
we are looking for curves which produce $p$-planes in $\cc^{m+p}$
which meet given $m$-planes sampled at prescribed values for~$s$.
In Figure~\ref{figdomains} we show the transition between the time
and the frequency domain.

\begin{figure}[hbt]
\begin{center}
\begin{picture}(240,100)(0,0)

\put(10,70)
{
  $\left\{
      \begin{array}{l}
         \dot{{\bf x}} = A {\bf x} + B {\bf u} \\
         {\bf y} = C {\bf x}
      \end{array}
    \right.$
}
\put(0,0)
{
  $\left[
      \begin{array}{c}
          C(sI_n - A)^{-1}B \\
              I_m
      \end{array}
   \right]$
}

\put(140,70)
{
  $\left\{
      \begin{array}{l}
         \dot{{\bf z}} = F {\bf z} + G {\bf y} \\
         {\bf u} = H {\bf z} + K {\bf y}
      \end{array}
   \right.$
}

\put(140,0)
{
  $\left[
      \begin{array}{c}
              I_p \\
          H(sI_q - F)^{-1}G + K \\
      \end{array}
   \right]$
}

\put(50,50){\vector(0,-1){25}}
\put(53,35){(1)}

\put(180,25){\vector(0,1){25}}
\put(183,35){(3)}

\put(100,0){\vector(1,0){30}}
\put(110,3){(2)}

\end{picture}
\end{center}
\caption{Transitions between time and frequency domain,
         as defined in Table~\ref{tabdomains}.}
\label{figdomains}
\end{figure}

When controlling a machine with $n$ internal states
with a controller using $q$ internal states, we can
place $n+q$ poles.  On the other hand, the dimension
of the geometric problem is $m p + q (m+p)$.
Depending on values for $n$, $m$, $p$, and $q$,
we distinguish three cases:

\begin{description}

  \item[1) $n+q < m p + q (m+p)$ underdetermined:]

For a generic machine, there is a set of feedback laws.
The set has dimension $m p + q (m+p) - n - q$, 
and has degree $d(m,p,q)$, i.e: for a generic choice of
the parameters, we have $d(m,p,q)$ complex feedback laws.
This degree can easily be computed combinatorially.

  \item[2) $n+q = m p + q (m+p)$ dimension zero:]

For a generic machine, there are exactly $d(m,p,q)$
complex feedback laws.  Every feedback law places
all $n+q$ poles at the desired locations.  It may be
that no feedback law has all its coefficients real,
see~\cite{EG02a,EG02b}.

  \item[3) $n+q > m p + q (m+p)$ overdetermined:]

For a generic machine, there are no feedback laws
which place all $n+q$ poles at the desired locations.

\end{description}

The numbers $m$, $p$, $n$ are fixed, given on input.
We can choose $q$ to arrive always in a favorable condition and
get feedback laws. In particular, 
if a system is overdetermined for static output feedback ($n>mp$),
we can always find $q$ to make it convert to underdetermined or 
dimension zero case with dynamic feedback. This is because $q(m+p)$ 
always 
increases faster than $q$.

Presently, we resolve the undetermined case by choosing
additional input planes to the geometric problem.
For nongeneric machines, the dimension of the set of
feedback laws is higher than expected.  Recent advances
with homotopies (see e.g. ~\cite{SV,SVW1,SVW4,SW})
allow to treat positive dimensional solution sets, 
but we defer the application of these
recent homotopy methods to a future paper.

\section{Symbolic-Numeric Calculations}

To transform a solution obtained by the homotopy methods into
a description of a machine to execute the dynamic feedback law,
we need to calculate the Smith normal form to compute the inverse
of a matrix with polynomial entries.
More precisely, the output of the homotopies as an $(m+p)$-by-$p$ matrix of
polynomials in $s$:
\begin{equation}
   \left[
      \begin{array}{c}
          U(s) \\ V(s)
      \end{array}
   \right], 
\end{equation}
where $U(s)$ is a $p$-by-$p$ matrix and $V(s)$ is an $m$-by-$p$ matrix
of polynomials in $s$, satisfies
\begin{equation}
  \det\left[
      \begin{array}{cc}
         U(s)   &     C(sI_n - A)^{-1}B \\
         V(s)   &        I_m  \\
      \end{array}
    \right] = 0, \label{eqhomotopyout}
\end{equation}
for the given poles.
We can right multiply~(\ref{eqhomotopyout}) by
\begin{equation}
   \left[
      \begin{array}{cc}
         U^{-1}(s) & 0 \\
           0       & I_m \\
      \end{array}
   \right].     \label{eqmultiplier}
\end{equation}
The result of this multiplication is
\begin{equation}
   \det\left[
      \begin{array}{cc}
              I_p         &     C(sI_n - A)^{-1}B \\
          V(s) U^{-1}(s)  &        I_m            \\
      \end{array}
   \right] = 0.      \label{eqmultiplyresult}
\end{equation}
Since the multiplier matrix~(\ref{eqmultiplier}) is of full rank,
its determinant is nonzero and the original intersection condition 
remains.
This multiplication does not affect the input conditions,
which are at the right part of ~(\ref{eqgeometric}).
By comparing~(\ref{eqgeometric}) with~(\ref{eqmultiplyresult}),
we can apply the realization algorithms to extract 
$(F,G,H,K)$ from the matrix $V(s) U^{-1}(s)$. 

In the next subsection we show how the calculation of a Smith
normal form requires the calculation of greatest common divisor.  

\subsection{Numerical Smith Normal Form}

The Smith normal form of a matrix of polynomials can be
seen as the discrete analogue of the singular value decomposition.
More precisely, for any $n$-by-$m$ matrix $A(x)$ whose
entries are polynomials in~$x$, there exists a unimodular $n$-by-$n$
matrix $P(x)$ and a unimodular $m$-by-$m$ matrix $Q(x)$ so that
\begin{equation}
   P(x) A(x) Q(x) = D(x),
\end{equation}
where $D(x)$ is an $n$-by-$m$ matrix which has only nonzero 
polynomials in $x$ on its diagonal.  Furthermore, denoting the
$i$th element on the diagonal of $D(x)$ by $D_i$, we have that
$D_i$ divides $D_{i+1}$.  
The matrix $D(x)$ is called the Smith normal form of $A(x)$.
Since unimodular matrices are invertible, we can rewrite $A(x)$ as
\begin{equation}
   A(x) = P^{-1}(x) D(x) Q^{-1}(x),
\end{equation}
which reveals the following expression for the inverse of $A(x)$:
\begin{equation}
   A^{-1}(x) = Q(x) D^{-1}(x) P(x),
\end{equation}
which of course only exists if $D(x)$ has full rank. We get the inverse
$D^{-1}(x)$ by inverting every entry on the diagonal of $D(x)$.

The Smith normal form can be computed by solving the extended 
greatest common divisor problem.
In particular, we wish to find polynomials $k(x)$ and $l(x)$
so that the following equality holds:
\begin{equation} \label{eqbezeq}
   d(x) = {\rm GCD}(a(x),b(x)) = k(x) a(x) + l(x) b(x). 
\end{equation} 
The calculation of the greatest common divisor is used to
reduce columns
\begin{equation}
  \left[
     \begin{array}{rr}
         k(x) & l(x) \\
         - \frac{b(x)}{d(x)} & \frac{a(x)}{d(x)}
     \end{array}
  \right]
  \left[
     \begin{array}{c}
        a(x) \\ b(x)
     \end{array}
  \right]
  =
  \left[
     \begin{array}{c}
        d(x) \\ 0
     \end{array}
  \right]
\end{equation}
or to reduce rows
\begin{equation}
  \left[
     \begin{array}{rr}
        a(x) & b(x)
     \end{array}
  \right]
  \left[
     \begin{array}{rr}
         k(x) & - \frac{b(x)}{d(x)} \\
         l(x) & \frac{a(x)}{d(x)}
     \end{array}
  \right]
  =
  \left[
     \begin{array}{rr}
        d(x) & 0
     \end{array}
  \right].
\end{equation}
The matrices we use in the column and row reductions have determinant
equal to one, they are called unimodular.

Collecting the column reductions in $P(x)$ and the row reductions
in $Q(x)$, we can reduce any polynomial matrix $A(x)$ to a diagonal
form~$D(x)$, and so obtain the Smith normal form of $A(x)$.

\subsection{Numerical Greatest Common Divisor}

Suppose we are given two polynomials in one variable
and with complex coefficients, only known approximately.
We denote the given polynomials by $a(x)$ and $b(x)$,
and $d(x) = {\rm GCD}(a(x),b(x))$, with $\deg(d(x)) = r$.

If we apply the algorithm taught in elementary school to compute
the greatest common divisor of two natural numbers, then we
repeatedly have to divide polynomials.  This repetitive division
is numerically unstable as it involves the subtraction of
polynomials, which may lead to a dramatic loss of significant
digits when the polynomials have coefficients of equal magnitude.
Our approach to computing the greatest common divisor of polynomials
with approximate coefficients hinges on two operations: root finding
and interpolation.  Any introductory course in numerical analysis
describes numerically stable algorithms for these two operations.

For a given tolerance $\epsilon > 0$,
we define the numerical greatest common divisor of $a(x)$ and $b(x)$
as the monic polynomial whose roots are common to $a(x)$ and $b(x)$
within the given tolerance~$\epsilon$.  More precisely, if
$a(\alpha_i) = 0$, for $i=1,2,\ldots,\deg(a(x))$ and $b(\beta_i) = 0$,
for $i=1,2,\ldots,\deg(b(x))$, then we can rearrange the indices of
the roots of the two polynomials so that the $r$ common roots
appear with lowest indices.  Then we write

\begin{equation}
   a(x) = \prod_{i=1}^r (x - \alpha_i) 
          \prod_{i=r+1}^{\deg(a(x))} (x - \alpha_i)
        = d_1(x) \prod_{i=r+1}^{\deg(a(x))} (x - \alpha_i)
\end{equation}
and
\begin{equation}
   b(x) = \prod_{i=1}^r (x - \beta_i)
          \prod_{i=r+1}^{\deg(b(x))} (x - \beta_i)
        = d_2(x) \prod_{i=r+1}^{\deg(b(x))} (x - \beta_i)
\end{equation}
where $|\alpha_i - \beta_i| \leq \epsilon$, for $i=1,2,\ldots,r$,
and $|\alpha_i - \beta_j| > \epsilon$, for all $i$ and $j$ with
index higher than $r$.  The polynomials $d_1(x)$ and $d_2(x)$
are numerical approximations for the greatest common divisor
$d(x)$ of $a(x)$ and $b(x)$.

Now that we have numerical approximations for $d(x)$,
we want to find $k(x)$ and $l(x)$ defined in~(\ref{eqbezeq}).
We determine $k(x)$ by interpolation
at those roots of $b(x)$ not shared by $a(x)$ replacing 
in~(\ref{eqbezeq}) $x$ by $\beta_i$, for $i= r+1,\ldots,\deg(b(x))$:
\begin{equation} \label{eqintcondk}
   d(\beta_i) = k(\beta_i) a(\beta_i)
   \quad {\rm or} \quad
   k(\beta_i) = \frac{d(\beta_i)}{a(\beta_i)},
   \quad {\rm for}~i > r.
\end{equation}
Note that as $i > r$: $a(\beta_i) \not= 0$.
The interpolation conditions in~(\ref{eqintcondk}) determine
$k(x)$ uniquely as a polynomial of degree $\deg(b(x))-r-1$.
Analogously,
we determine $l(x)$ by interpolation at those roots of $a(x)$
not shared by $b(x)$ replacing in the equality~(\ref{eqbezeq})
$x$ by $\alpha_i$, for $i= r+1,\ldots,\deg(a(x))$:
\begin{equation} \label{eqintcondl}
   d(\alpha_i) = l(\alpha_i) b(\alpha_i)
   \quad {\rm or} \quad
   l(\alpha_i) = \frac{d(\alpha_i)}{b(\alpha_i)},
   \quad {\rm for}~i > r.
\end{equation}
Note that as $i > r$: $b(\alpha_i) \not= 0$.
The interpolation conditions in~(\ref{eqintcondl}) determine
$l(x)$ uniquely as a polynomial of degree $\deg(a(x))-r-1$.

\subsection{Numerical Experiments}

The algorithms described above have been implemented in C.
In this section we list some numerical results, obtained
by practical comparisons between our new algorithm and the
elementary approach, for random and specific input data.

We call the elementary school algorithm by repetitive
division the ``naive algorithm'' and the algorithm with
root finding the ``advanced algorithm''.
When the tolerance of $10^{-8}$ is used to decide whether two numbers 
are equal, the naive algorithm runs much faster than the other one,
while the advanced algorithm is more numerically stable
when the polynomial degree is less than 30 and the degree of the GCD 
is less than~15. Although the advanced algorithm is relatively slow, 
the time spent for each GCD computation is just trivial,
less than~10 milliseconds even for the polynomial of degree~30
and the GCD of degree~15.
Table~\ref{tabgcd1} summarizes this experiment.

\begin{table}[hbt]
\begin{center}
\begin{tabular}{|c|c|c|c|} \hline
 degree of $a$ and $b$  & degree of GCD & naive$(\%)$ & advanced$(\%)$ 
\\
 \hline \hline
        ~5          &      ~3       &      100.0   &    100.0  \\
        10          &      ~5       &      ~99.2   &    100.0  \\
        15          &      ~8       &      ~99.5   &    100.0  \\
        20          &      10       &      ~99.5   &    100.0  \\
        25          &      13       &      ~98.2   &    100.0  \\
        30          &      15       &      ~88.2   &    100.0  \\ 
\hline
\end{tabular}
\caption{Comparison of the naive algorithm and the advanced 
         algorithm on random data.
The first column lists the degree of the given polynomials $a$ and~$b$.
The second column lists the degree of the GCD.
The data in columns three and four shows the percentage of success 
for 1000 tests of the two algorithms.}
\label{tabgcd1}
\end{center}
\end{table}

We also compared the performance of the algorithms for some specific 
data.

In the first experiment, the leading coefficient of the polynomial~$b$
is $10^{-5}$ times leading coefficient of the polynomial~$a$.
We found that the naive algorithm fails completely 
when the degree~$e$ of the input polynomials and the
degree~$r$ of the GCD satisfy the following relation:
\begin{equation}
   e = \left\{
         \begin{array}{ll}
            2 r - 1 & \mbox{for } r \mbox{ is odd} \\
            2 r     & \mbox{for } r \mbox{ is even}
         \end{array} 
       \right.
\end{equation}
Other relations between the degrees work fine.
While the advanced algorithm remains numerically stable 
no matter what the relations between the degrees are.
See the column with header specific numbers(1) in 
Table~\ref{tabgcd2} for the experiment of this case. 

Second, we test with some other specific numbers, say,
if the higher degree coefficients of two polynomials are 
very near each other, then the two algorithms perform totally different.
In particular, we tested polynomials with equal leading coefficients
and with difference between the second-highest coefficient equal
to~$10^{-5}$.
When the degree of the input polynomials
is smaller than~5 and the degree of the GCD is smaller than~3,
both algorithms work perfectly. 
When the degree of the input polynomials is larger than~10 
and the degree of GCD is larger than~5, the naive algorithm 
fails completely, but the advanced algorithm shows the same
numerical stability as on the random number case.  
See columns with header specific numbers(2) 
in Table~\ref{tabgcd2} for the experiment of this case. 

\begin{table}[hbt]
\begin{center}

\begin{tabular}{|c|c|c|c|c|c|} \hline
   \multicolumn{2}{|c|}{~}
            & \multicolumn{2}{c|}{specific numbers$(1)$} 
            & \multicolumn{2}{c|}{specific numbers$(2)$} \\ \hline
   degree   & degree & naive  & advanced & naive  & advanced \\
 of $a$ and $b$ & of GCD & method$(\%)$ & method$(\%)$ & method$(\%)$ 
            & method$(\%)$ \\ \hline \hline
   ~5 & ~3 & ~0.0 & 100.0 & 99.9 & 100.0 \\
   10 & ~5 & 99.1 & 100.0 & ~0.1 & 100.0 \\
   15 & ~8 & 98.8 & 100.0 & ~0.0 & 100.0 \\
   20 & 10 & ~0.3 & 100.0 & ~0.0 & 100.0 \\ \hline
\end{tabular}
\caption{Comparison of the naive algorithm and the 
         advanced algorithm on specific data. 
The first column lists the degree of the given polynomials $a$ and~$b$.
The second column contains the degree of the GCD. 
The data in columns three to six shows the percentage of success
for 1000 tests of the two algorithms. }
\label{tabgcd2}
\end{center}
\end{table}

These experiments provide us with practical evidence that
for both random and specific inputs,
the advanced algorithm shows its strong numerical stability. 
Concerning the speed of the GCD algorithm, the time needed of the 
realization algorithm is negligible compared to the calculation time
of the feedback laws with homotopies.

\section{Realization of Multi-Input Multi-Output Systems}

In the previous section, we gave the derivation of the transfer function
of the dynamic compensator as the output of the homotopies.
We will give a modified algorithm based on~\cite[pages 389--416]{AN97}
to obtain minimal 
or irreducible realizations, which realize a system with the least number 
of dynamic elements. These modifications were made to fit
the output format of the software used to compute the feedback laws. The
necessity of the modifications will be discussed at the end of this section.
We will show how to obtain realizations 
$\{F_c, G_c, H_c, K_c\}$ of the transfer function $T(s)$ 
in controller form first.
Then, we will use the property of the output of
homotopies to show the realizations are irreducible, so they are also 
controllable.

From~(\ref{eqmultiplyresult}), the transfer function can be written as
\begin{equation}
       T(s) = V(s)U^{-1}(s).
\end{equation}
In accordance with convention, we would replace $V(s)$ by $N(s)$ which 
stands for numerator, replace $U(s)$ by $D(s)$ which stands for 
denominator.

According to Theorem~3.3 in~\cite[page 391]{AN97},
realizations exist if and only if
$T(s)$ is a matrix of rational functions and satisfies 
\begin{equation}
  \lim_{s\rightarrow\infty}T(s) < \infty, 
\end{equation}
i.e., if and only if $T(s)$ is a proper rational matrix.
Given the transfer function matrix $T(s) = N(s)D^{-1}(s)$ as a $(m\times p)$
 proper rational matrix. Let $d_j$ = the highest degree of $j$th column in 
the $D(s)$ ($d_j\geq 0, j=1,2,\ldots, p$). Define 
\begin{equation}
          \Lambda(s) = \mbox{diag}(s^{d_1},\ldots,s^{d_p}),  
\end{equation}
and
\begin{equation}
          S(s) =\mbox{block diag}
                        \left(
                         \begin{array}{l}
                           \left[
                             \begin{array}{c}
                                1\\
                                s\\
                                \vdots\\
                                s^{d_j-1}
                             \end{array}
                           \right]
                           j = 1,\ldots, p
                         \end{array}
                       \right).
\end{equation}
If $d_j=0$, just skip that column and continue to fill the next column 
of the $S(s)$ matrix.
Note that $S(s)$ is an $q(=\sum_{j=1}^p d_j)\times p$ polynomial 
matrix.
Write
\begin{equation} \label{eqfac}
   D(s) = D_h\Lambda(s) + D_lS(s)
\end{equation}
$D_h$ is the highest column degree coefficient matrix of $D(s)$.
For example, if $D(s) =\left[
         \begin{array}{l}
           \begin{array}{c}
              3s^2+1\\
              2s
           \end{array}
           \begin{array}{c}
              2s\\
              s  
           \end{array}
         \end{array}   
        \right],$
then the highest column degree coefficient matrix 
$D_h =\left[
         \begin{array}{l}
           \begin{array}{c}
              3\\
              0
           \end{array}
           \begin{array}{c}
              2\\
              1  
           \end{array}
         \end{array}   
        \right],$ 
and $D_lS(s)$ given in~(\ref{eqfac}) accounts for the remaining lower 
column 
degree terms $D(s)$, with $D_l$ being a matrix of coefficients.

  In general, $|D_h|\neq0$, and define $p\times p$ and $p\times q$ 
matrices
\begin{equation}
    G_p = D_h^{-1}, \hspace{10mm}   F_p = -D_h^{-1}D_l,
\end{equation}    
respectively. Then $F_c,~G_c$ can be determined from
\begin{equation}
   F_c = \overline{F}_c + \overline{G}_cF_p,\hspace{10mm} 
   G_c = \overline{G}_cG_p, 
\end{equation}
where $\overline{F}_c =\mbox{block diag}[F_1, F_2 \ldots, F_p]$ with

\begin{equation}
   F_j =\left[ 
          \begin{array}{cccc}
              0   &   &               &   \\
           \vdots &   & I_{d_{j-1}}   &   \\
              0   & 0 &  \ldots       & 0 \\
          \end{array}
        \right]
       \in R^{d_j\times d_j}, 
\end{equation}
\begin{equation}
   \overline{G}_c = \mbox{block diag}
                        \left(
                           \left[
                              \begin{array}{c}
                                0\\
                                \vdots\\
                                0\\
                                1
                              \end{array} 
                           \right]
                           \in R^{d_j}, j = 1, \ldots, p
                        \right).     
\end{equation}

\noindent When $d_j =0$, we just skip the corresponding $F_j$ and 
continue 
to fill the $\overline{F}_c$ matrix with the $F_{j+1}$ matrix; 
we also need to add 
a zero column at the $j$th column of the $\overline{G}_c$ matrix.

Then we can determine $H_c$ and $K_c$ such that
\begin{equation} \label{eqNs}
   N(s) = H_cS(s)+K_cD(s),
\end{equation}
and note that
\begin{equation}
   K_c = \lim_{s\rightarrow \infty}T(s).
\end{equation}
Therefore, $H_c$ can be determined from~(\ref{eqNs}).

An $q$th-order realization of $T(s)$ in controller form
is now given by the equations
\begin{equation}
  \dot{\bz}_c = F_c \bz_c + G_c \by, \quad 
   {\bf u} = H_c \bz_c+K_c \by.
\end{equation}

According to the format of the output of the software, 
$q = \sum_{j=1}^p d_j$ is equal to the 
minimal order of the dynamic compensator. 
Therefore, this algorithm gives us a minimal realization of 
the transfer function matrix $T(s)$ and the result is also observable.

The main difference between the modified algorithm and the original algorithm
given in~\cite{AN97} is that the original algorithm limits  $d_j \geq 1$,
 while our modified algorithm works for $d_j \geq 0$, where $d_j$ is the
 highest column degree of $j$th column in the $D(s)$. Some $d_j$ must be 
equal to zero when the number of output is larger than the order 
of the dynamic compensator. In this case the modified algorithm become 
necessary. The correctness of the modified algorithm is 
verified with experiments.   

\section{Software}

The dynamic feedback laws were calculated with the aid
of PHCpack~\cite{Ver99}.
While the second public release of PHCpack implemented 
the dynamic pole placement problem in its geometric form,
additional software had to be written, concerning:
\begin{enumerate}
   \item[0.] a limit on the number of feedback laws;
   \item[1.] an interface between Ada and C; and
   \item[2.] a collection of C routines for the realization.
\end{enumerate} 
The limit on the number of feedback laws was imposed as a matter
of convenience, to control the practical complexity.
We elaborate the other two items in the following subsections.

\subsection{A C interface to PHCpack}

PHCpack is written in Ada, while the programs to process the
feedback laws are in the lower level language~C.

We can build a portable interface to the Ada routines in
PHCpack with C functions because the language Ada
has the {\tt pragma Import} construction to 
call routines from other languages such as C and it supports
conversions for C integers, doubles, and arrays of these C types.
Furthermore, the gnu-ada compiler provides a mechanism to call 
Ada routines from a C main program and to call C functions
from Ada.  As the gnu-ada compiler is integrated in the gcc 
compilation system, our interface is portable.
In particular, we ran our implementation successfully on SUN
workstations running Solaris and on PCs running Linux and Windows.

To exchange data efficiently, programs in Ada or C should
define exchange protocols of structured data types into 
basic data types for which automatic conversions are supported.
More precisely, we represent structured data types into arrays
of doubles and arrays of integers.
The language C is restricted in returning dynamically allocated
variables.  Therefore, data allocated in a C function is
passed by the C function calling an Ada function for further
processing of the data.

A typical sequence of calls goes as follows.  First a C function
gathers problem data and prepares the input to an Ada routine
of PHCpack.  The Ada routine, called from C, uses path tracking
to solve the problem, and then calls a C function to process
the results obtained with PHCpack.
So the C programmer who uses PHCpack should thus provide two
C functions: one to prepare the input and one to process the output.
This ``hand-in-glove'' interface is appropriate for a C programmer
collaborating with an Ada programmer (which is the case of the
authors), who only have to agree on the prototypes of the routines.

\subsection{The Organization of the Software}

\begin{figure}[hbt]
\centerline{\psfig{figure=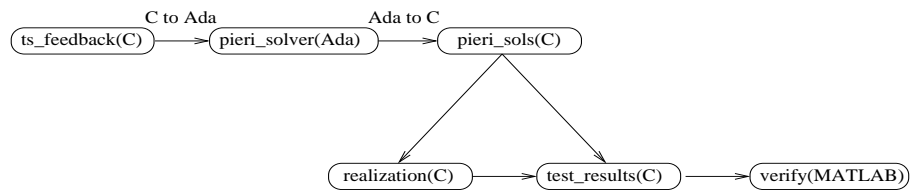,height=2.4cm,width=12cm}}
\caption{Organization of the Software}
\label{figorgsoft}
\end{figure}

In Figure~\ref{figorgsoft}, the arrows indicate the order
of function calls in the computation and realization of 
the dynamic feedback laws.  In this section we give a short
description for each of the procedures, some written in C,
others in Ada.

\begin{description}
  \item[ts\_feedback(C):]  ts\_feedback reads all the 
     input information from a file, including the number of the 
     internal states $n$, the input dimension $m$, the output 
     dimension $p$, the number of the internal states for the
     dynamic compensator $q$ and the number of output feedback.
     Also the user should give the $A, B, C$ matrices (or let the matrices
     be generated randomly) of the given plant and $n+q$ eigenvalues.
     ts\_feedback.c computes $C(sI_n - A)^{-1}B$ at the interpolation
     points as the input planes. 
     With C to Ada interface, we pass the arrays of the input planes and the 
     interpolation points to the pieri\_solver 
     (an Ada procedure in PHCpack).
  \item[pieri\_solver(Ada):] The pieri\_solver calculates the 
     corresponding dynamic output feedbacks and passes them to
     the C program pieri\_sols.c.
  \item[pieri\_sols(C):] With Ada to C interface, the arrays in Ada 
     form are converted to the form in C. Then pieri\_sols.c calls 
     realization.c and tests the results. 
  \item[realization(C):] We use modified realization algorithm based 
     on~\cite{AN97} to get the realization of the dynamic output feedback,
     organized as follows:
   \begin{enumerate}
        \item[(a)] Get the transfer function $T(s)=N(s)D(s)^{-1}$ from the 
      output of the Ada program. The inverse of a polynomial matrix is a 
      rational polynomial matrix and it is mainly done by Poly\_Smith.
      See (10) in section 3 for more detail about how Smith
      form can work for the inverse of a polynomial matrix.
        \item[(b)] The realization function implements the modified 
      realization algorithm to get a minimal realization of the dynamic 
      compensator.      
        \item[(c)] Evaluate the transfer function $T(s)$ at some random
      point and compare it with the result after realization 
      ($H(sI_q-F)^{-1}G+K$) at the same point.
      If the values are the same, the realization is correct. 
   \end{enumerate}

\item[test\_results(C):] We can evaluate equation~(\ref{eqgeometric})
  at the given poles and calculate the determinant with the previous 
  result.  If the determinant is zero, the pole is the eigenvalue of the 
  closed-loop system.  As mentioned above, equation~(\ref{eqgeometric}) is 
  algebraically equal to equation~(\ref{eqrewrite}), which is 
  the characteristic equation of the closed-loop system.

\item[verify(MATLAB):]  Finally, a MATLAB script verifies
  the results by comparing the computed poles with the given poles
  and finding the condition number for each given pole.

\end{description}

\subsection{Availability of the Software}

The C routines for the numerical realization algorithms
written by Yusong Wang are available 
at {\tt http://www.math.uic.edu/{\~{}}jan}
in the distribution of release~2.2 of the 
source code of PHCpack (see~\cite{Ver99}, the first
version of this package is archived by Netlib).
In particular, the collection of routines can be found
in the directory ``Feedback'' of the source code.
Also the input data files for the applications discussed in
the next section are available in this directory.

\section{Applications}

In this section we illustrate the usefulness of our approach.
Some application examples will show it is necessary to 
find some dynamic output feedback laws for some specific situations,
especially for the cases there is no feedback laws or no real feedback
laws at some given poles for the static feedback problem. 

We assume that the input data for our applications is given
as $(A,B,C)$, i.e.: a triple of three matrices of the linear
system $\dot{\bx} = A \bx + B \bu$, and $\by = C \bx$, where
$\bx, \bu$, and $\by$
are vectors of internal states, input, and output respectively.

\subsection{Satellite Trajectory Control}

This application concerns the design of output feedback laws
to keep a satellite in orbit.  We treated this problem in our
previous paper~\cite{VW02} with static output feedback and with
ad hoc MATLAB scripts to aid the manipulations.
In~\cite{VW02}, we already found two real output feedback laws for 
static case.
 Now we want to get some feedback laws for dynamic case as an extension 
of our
 method. 

The linearized state-space equations for the satellite problem are 
given below:
\begin{equation}
  \begin{array}{l} 
         A = \left[
                \begin{array}{cccc}
                    0      &  1      &  0    & 0       \\
                    0.3578 &  0      &  0    & 0.8525  \\
                    0      &  0      &  0    & 1       \\
                    0      & -0.5596 &  0    & 0       \\   
                \end{array}
             \right];\\
                \\     
         B = \left[
                \begin{array}{cc}
                    0      &  0      \\
                    1.3411 &  0      \\
                    0      &  0      \\
                    0      &  1.0867 \\
                \end{array}
             \right].\\
                \\     
 \end{array}           
\end{equation}
We choose value for each parameter randomly. While this may seem 
unrealistic, our choice can be
justified by an appropriate selection of units.                      
More detail about the model can be found in~\cite{DB98,Kai80}.
We define $C$ as some random matrix -- which can be interpreted as a 
random
 projection of the states onto a plane.

For this satellite example, we have $n=4$, $m=p=2$,  which is dimension 
zero for the
static case. When $q=1$, the system become underdetermined since $n+q=5 
< mp+q(m+p)=8$ and there are three degrees of freedom. We choose the 
eigenvalues as
$(\frac{-2+i}{\sqrt{5}}, \frac{-2-i}{\sqrt{5}}, -5, -7, -3.0, 
-0.1068,-0.7834,
-0.9582)$. The last three of eigenvalues are randomly selected. We 
found two real feedback laws and six 
complex complex feedback laws.
 We also use the MATLAB script to verify the results and calculate the 
condition number for the eigenvalue~$\lambda_i$ by
\begin{equation}\label{eqcond}
   \frac{1}{|\by_i^H \bx_i|},
   \quad i=1,2,\ldots,n+q,
\end{equation}  
where vectors $\bx_i$ and $\by_i$ denote the unit right and left 
eigenvectors of the closed-loop system. See~\cite[page 323]{GV96} for the 
derivation of $|\by_i^H \bx_i|$ as the reciprocal of the condition number 
for the eigenvalue $\lambda_i$.
 By substituting the result into the closed-loop system, we find the relative 
difference of the computed eigenvalues vs. the given eigenvalues is bounded 
by $10^{-11}$ and the order of condition numbers are at most $10^3$, calculated
 with~(\ref{eqcond}).
 The total CPU time spent on the calculation of the dynamic feedback 
laws is 2 seconds and 320 milliseconds on a 2.4GHz workstation running 
Linux System.

\subsection{Numerical Example}
In this section we report on two numerical examples in
the literature~\cite{RW96} \cite{Wan96}.
\begin{description}
  \item[Numeric Example A:]
We will use the example in~\cite{RW96} 
to illustrate the following situation: 
when a system is overdetermined for static output 
feedback ($n > mp$, $q=0$), for which no feedback laws can be found at the
desired poles, we can choose a q to make
it convert to underdetermined or dimension zero case 
($n+q \leq mp + q(m+p)$), for which we can find feedback laws.

Consider the system, given by
\begin{equation}
  \begin{array}{l} 
     A = \left[
                \begin{array}{rrrrrr}
                   0 & 0 & 0  & 0 & 0 & 0  \\
                   1 & 0 & 0  & 0 & 0 & -1 \\
                   0 & 1 & 0  & 0 & 0 & 0  \\
                   0 & 0 & 0  & 0 & 0 & 0  \\
                   0 & 0 & 0  & 1 & 0 & 0  \\
                   0 & 0 & -1 & 0 & 1 & 0  \\
                \end{array}
             \right] ; \\ \\
     B = \left[
                \begin{array}{rr}
                  1 &  3 \\
                  0 &  0 \\
                  0 & -1 \\
                  0 &  1 \\
                  0 &  1 \\
                  0 &  0 \\
                \end{array}
             \right] ; \quad
         C = \left[
                 \begin{array}{cccccc}
                   0 & 0 & 1 & 0 & 0 & 0 \\  
                   0 & 0 & 0 & 0 & 0 & 1 \\
                 \end{array}
             \right].
 \end{array}           
\end{equation}
Here, $m=p=2$, $n=6$, therefore $n > mp$ which is overdetermined
for static case.  We will choose a $q$ to make the system 
underdetermined.  The minimum possible $q$ is 1, 
then we have $n+q=7 < mp+q(m+p)=8$, so there is one degree of freedom.
With the calculation of PHCpack and the C interface, we can easily
get some dynamic compensators of McMillan degree $q=1$ to control 
the system. We verified the result by comparing the given and 
computed eigenvalues with a MATLAB script.
For a choice $(-0.1,-1.5,-0.9,-0.7,-6.0,-3.5,-8.0)$ of 7 eigenvalues 
and an additional pole $-0.1053$ which is generated by random, 
the relative difference of the computed eigenvalues vs. the given eigenvalues
of the closed-loop system is bounded by $10^{-9}$.  In this specific example,
we find 8 solutions and 4 of them are real (sometimes 6 are real,
depending on the additional input plane which are randomly generated).
The order of the condition number computed with~(\ref{eqcond}) is 
no more than $10^3$ for all of the given eigenvalues.
The total CPU time spent for this numeric example is around 3 seconds 
with the workstation mentioned above.

  \item[Numeric Example B:]
The second example can be found in~\cite[Example 3.7]{Wan96}, 
with $n=8$ and $m=p=3$.

The system is given as below
\begin{equation}
  \begin{array}{l} 
         A = \left[
                \begin{array}{rrrrrrrr}
                      0  & -1  &  0  &  0  &  0  &  0  &  0  &  1 \\
                      1  &  2  &  0  &  0  &  1  &  0  &  0  & -2 \\
                      0  & -1  &  0  &  0  &  5  &  0  &  0  &  0 \\
                      0  &  0  &  1  &  0  & -7  &  0  &  0  & -2 \\
                      0  & -1  &  0  &  1  &  4  &  0  &  0  &  2 \\
                      0  & -2  &  0  &  0  &  2  &  0  &  0  &  3 \\
                      0  &  0  &  0  &  0  & -1  &  1  &  0  & -2 \\
                      0  & -1  &  0  &  0  &  1  &  0  &  1  & -1 \\
                \end{array}
             \right];\\
                \\     
         B = \left[
                \begin{array}{rrr}
                      0  &  1  &  2 \\
                      1  &  0  &  1 \\
                     -1  & -1  & -3 \\
                      1  &  0  &  1 \\
                      0  &  2  &  4 \\
                      2  &  1  &  5 \\
                     -1  &  1  &  1 \\
                      1  & -1  & -1 \\
                \end{array}
             \right]; \\ \\
         C = \left[
                \begin{array}{rrrrrrrr}

                      0  &  1  &  0  &  0  & -2  &  0  &  0  &  0 \\
                      0  &  0  &  0  &  0  &  1  &  0  &  0  &  0 \\
                      0  &  0  &  0  &  0  &  0  &  0  &  0  &  1 \\
                \end{array} 
             \right].\\
 \end{array}           
\end{equation}
This system is underdetermined for static output feedback ($n < mp, q=0 
$).
When the given poles are $( -0.8090 + 0.5878i,  -0.9511 - 0.3090i, 
-0.3090 - 0.9511i,  -0.3090 + 0.9511i,  -0.9511 + 0.3090i,  
-0.8090 - 0.5878i,  -0.5878 + 0.8090i, -0.5878 - 0.8090i,
-0.1883)$, in which the first 8 are picked at the unit circle 
and the last one is a negative real number chosen by random,
we find 42 feedback laws with our software and 4 of them are real.
The relative difference of the computed eigenvalues vs. the given
eigenvalues is $10^{-9}$.  The order of the condition number is 
bounded by $10^5$ if computed with~(\ref{eqcond}).
The total CPU time spent on this example is around 50 seconds 
to find all 42 feedback laws.  It only takes 1 second and 180 
milliseconds if user just needs one feedback law.

\end{description}

\subsection{Aircraft Control}
It may be that for a given selection of poles, all static feedback laws 
have
coefficients with nonzero imaginary parts.  In this case, 
we will design dynamic output feedback laws, and exploit the
additional freedom to place the poles at their originally
selected locations.  The model from MathWorks~\cite{MW01}
could be used to illustrate this kind of situation.
It also can be found in any standard text in aviation 
for a more complete discussion of the physics behind aircraft 
flight. The state-space equations are given below.

The jet model during cruise flight at MACH=0.8 and H=40,000ft. is
\begin{equation}
  \begin{array}{l} 
         A = \left[
                \begin{array}{cccc}
                   -0.0558  & -0.9968  & ~~0.0802  & 0.0415  \\
                   ~~0.5980 & -0.1150  & -0.0318   & 0       \\
                   -3.0500  & ~~0.3880 & -0.4650   & 0       \\
                    ~0      & ~~0.0850 &  ~~1.0000 & 0       \\   
                \end{array}
             \right];\\
                \\     
         B = \left[
                \begin{array}{cc}
                    ~~0.0073 &  0      \\
                   -0.4750   &  0.0077 \\
                    ~~0.1530 &  0.1430 \\
                    ~0       &  0      \\
                \end{array}
             \right]; \quad
         C = \left[
                \begin{array}{cccc}
                    0  &  1  &  0  &  0 \\
                    0  &  0  &  0  &  1 \\
                \end{array} 
             \right].\\

 \end{array}           
\end{equation} 
For this example, $n=4$, $m=p=2$, so $n=mp$ when $q=0$, 
which is the dimension zero case.  We can find some feedback laws
with PHCpack.  When the chosen poles are 
$(-0.234,-1+3.2i,-1-3.2i,-3.0)$,
we found two complex feedback laws.  The total time spent 
is 140 milliseconds and the realization part just needs trivial time.
By computing the eigenvalues of the closed-loop system with 
the MATLAB script, we find the difference between the computed
eigenvalues and the given eigenvalues is just $10^{-14}$ and the
condition number is around~10.  Although the result seems quite nice,
we still want to find some real feedback laws.  We could choose $q=1$
to find some dynamic feedback laws, therefore $n+q=5 < mp+q(m+p)=8$,
then we have three degrees of freedom.  With the calculation of 
PHCpack and C interface, we get 8 solutions, including 2 or 4, or 
6 real solutions, depending on different additional eigenvalues. 
For a choice of $(-0.234,-1+3.2i,-1-3.2i,-3.0,-7.0)$, in which
the first four are the same as in the static case, and three randomly
generated poles are $(-0.944, -0.995, -0.904)$, the relative 
difference of the computed eigenvalues vs. the given eigenvalues is
 bounded by $10^{-10}$. The total CPU time is 2 seconds and 680 milliseconds.
For this choice of poles, we found four real feedback laws. 
For most of the 8 solutions, the condition numbers of the closed-loop
system are less than $10^2$, and few of the solutions have the
condition number $10^{4}$, calculated with formula~(\ref{eqcond}).

\section{Conclusions}

In this paper we showed the practical feasibility of computing
dynamic feedback laws using numerical homotopy algorithms,
and we applied our software to examples from the control literature.
Our publicly available implementation of the numerical greatest
common divisor and numerical Smith normal form may be of
independent interest.

\end{document}